\documentclass[sn-mathphys-num]{sn-jnl}

\usepackage{graphicx}%
\usepackage{multirow}%
\usepackage{amsmath,amssymb,amsfonts}%
\usepackage{amsthm}%
\usepackage{mathrsfs}%
\usepackage[title]{appendix}%
\usepackage{xcolor}%
\usepackage{textcomp}%
\usepackage{manyfoot}%
\usepackage{booktabs}%
\usepackage{algorithm}%
\usepackage{algorithmicx}%
\usepackage{algpseudocode}%
\usepackage{listings}%

\usepackage{lmodern}

\newtheorem{theorem}{Theorem}%

\newcommand{\clip}{\texttt{G}}
\newcommand{\mbR}{\mathbb{R}}
\newtheorem{example}{Example}%
\newtheorem{remark}{Remark}%

\newtheorem{definition}{Definition}%

\newcommand{\abstractText}{
In practical engineering and optimization, solving multi-objective optimization (MOO) problems typically involves scalarization methods that convert a multi-objective problem into a single-objective one. While effective, these methods often incur significant computational costs due to iterative calculations and are further complicated by the need for hyperparameter tuning.

In this paper, we introduce an extension of the concept of competitive solutions and propose the Scalarization With Competitiveness Method (SWCM) for multi-criteria problems. This method is highly interpretable and eliminates the need for hyperparameter tuning. Additionally, we offer a solution for cases where the objective functions are Lipschitz continuous and can only be computed once, termed Competitiveness Approximation on Lipschitz Functions (CAoLF). This approach is particularly useful when computational resources are limited or re-computation is not feasible.

Through computational experiments on the minimum-cost concurrent flow problem, we demonstrate the efficiency and scalability of the proposed method, underscoring its potential for addressing computational challenges in MOO across various applications.
}

\begin{document}

\title{$\gamma$-Competitiveness}
\subtitle{An Approach to Multi-Objective Optimization with High Computation Costs in Lipschitz Functions}

\author*[1,2]{\fnm{Ilgam} \sur{Latypov}}\email{i.latypov@iai.msu.ru}

\author[1,2]{\fnm{Dorn} \sur{Yuriy}}\email{dornyv@my.msu.ru}

\affil*[1]{\orgdiv{Institute for Artificial Intelligence}, \orgname{Lomonosov Moscow State University}, \orgaddress{\city{Moscow}, \country{Russia}}}
\affil[2]{\orgname{Moscow Institute of Physics and Technology}, \orgaddress{\city{Moscow},\country{Russia}}}

\abstract{\abstractText}

\maketitle

\section{Introduction}

In many practical applications, problems are often formulated with multiple criteria \cite{miettinen1999nonlinear,greco2006multiple,koksalan2011multiple}, meaning there is no multiple performance metrics that should be jointly optimized. To address such challenges, methods for multi-criteria decision-making and multi-objective optimization (MOO) have been developed \cite{koksalan2011multiple}. These approaches are widely used across various fields, including power supply and telecommunication networks \cite{altiparmak2006genetic,elmusrati2008applications,mastrocinque2013multi,bjornson2014multiobjective}, machine learning \cite{suttorp2006multi,zuluaga2013active,sener2018multi}, chemistry \cite{rangaiah2013multi}, biology \cite{boada2016multi}, and engineering \cite{marler2004survey}.

The multi-criteria optimization problem can be mathematically formulated as follows: 
\begin{align*}
    &\min_{x \in K} f \triangleq (f_1(x), ..., f_m(x))^T, & \tag{$T_0$}\label{opt:T0}\\
    &\text{s.t.} \quad x\in K.
\end{align*}

 where $K \subseteq \mathbb{R}^n$ is a non-empty, compact feasible set. In our analysis, we assume that the objective functions \( f_i: K \rightarrow \mathbb{R}^n_{++}\)  \((i = 1, \ldots, m)\) are Lipschitz continuous, and that the feasible set $K$ is convex. 

  A key challenge in solving multi-objective problems is determining how to minimize a vector-valued objective function. There are two common approaches to this. The first involves finding the set of Pareto-optimal solutions \cite{ngatchou2005pareto,konak2006multi, miettinen2008introduction, daulton2022multi}, which requires a detailed exploration of the feasible set. While this method is comprehensive, it is often computationally expensive and can produce redundant solutions. The second approach involves reducing the multi-objective problem to a single-objective one using an aggregation function—a process known as scalarization. Scalarization methods, such as weighting, target point, and direction-based methods, are well-documented in \cite{greco2006multiple}.
  
  However, scalarization typically relies on non-trainable parameters like weights, which must be predefined by experts. This introduces subjectivity and reduces the interpretability of the results. To overcome these limitations, we propose a new approach based on the concept of competitive solutions, inspired by classical algorithm analysis \cite{borodin2005online}. Our method avoids the need for weight-like parameters and provides a solution that is both interpretable and practical. 
 
 \begin{definition}[strictly $\gamma$-competitive solution]\label{def:strick}
    A point \( x \) is called a strictly \(\gamma\)-competitive solution for the set of functions \( f_i \) if \( \forall i = 1, \ldots, m \):
    $$
    f_i(x) \leq (1 + \gamma) \min_{y\in K} {f_i(y)}.
    $$
\end{definition}

\begin{definition}[relative $\gamma$-competitive solution]\label{def:gamma_competitive}
    Let there be a given set of feasible points $\{x_i\}_{i=1}^m$ $(x_i \in K, 1\leq i \leq m)$ with known values of corresponding metrics $\{f_i(x_i)\}_{i=1}^m$. Then a point $x$ is called a relative \(\gamma\)-competitive solution for the functions $\{f_i\}_{i=1}^m$ at $\{x_i\}_{i=1}^m$ if the following inequality holds:
    \begin{equation}
    f_i(x) \leq (1 + \gamma) f_i(x_i), \quad 1\leq i \leq m.
    \end{equation}
\end{definition}

Definition \ref{def:strick} is a natural extension of the concept of competitiveness. Specifically, a solution is called $\gamma$-competitive if, for each $1\leq i \leq m$ it satisfies $f_i(x) \leq (1+\gamma) f_i^*$. In other words it is $\gamma$-competitive solution for metric $f_i$ for each $i$. 

However, definition \ref{def:strick} has some limitations. The following example illustrates these limitations:

 \begin{example}\label{example:1}
    Suppose a company has operated over several periods using different strategies, with each strategy represented by a point \( x_i \). Metrics \( f_k(x_i) \) were calculated for these strategies. The company now needs to choose a strategy for the next period but lacks knowledge of future market conditions. It can, however, select a strategy that performs reasonably well based on past scenarios. 
\end{example}

In this case, we cannot apply definition \ref{def:strick} of competitiveness, as there is no information about the optimality of past actions. Therefore, we extend the definition to \ref{def:gamma_competitive}, which uses observable values as benchmarks.

In the remainder of this work, we focus on definition \ref{def:gamma_competitive}, as it generalizes definition \ref{def:strick}. When the optimal values $f_i(x_i^*) = f_i^*$ for each metric $i$ are known, we can use $x_i^*$ as reference points, making \ref{def:gamma_competitive} applicable.

Building on definition \ref{def:gamma_competitive}, we propose the Scalarization With Competitiveness Method (\ref{opt:T1}), which seeks to find a relative $\gamma$-competitive solution with the tightest feasible $\gamma$. We also introduce a computationally efficient method for approximately solving \ref{opt:T1} in cases where the objective functions are Lipschitz—referred to as the Competitiveness Approximation on Lipschitz Functions (\ref{opt:T3}). This method is particularly suited for scenarios when metric computation is computationally demanding or re-computation is not feasible.

This approach is especially practical for high-dimensional multi-criteria optimization problems with computationally intensive metrics. For instance, in telecommunication network optimization, one might aim to minimize both latency and bandwidth, which are derived from the minimum-cost flow and max-cut problems. Both metrics are objective values of linear optimization problems, which can be computationally demanding for large instances \cite{banos1995linear,martin2012large}.
 
\subsection{Paper organization}
In Section \ref{sec:task_statement} we formulate and discuss scalarization method, that based on the definition of relative $\gamma$-competitive solution -- \ref{opt:T1}. 
 In section \ref{sec:method} we introduce \ref{opt:T3} -- a method for finding an approximate solution to the scalarization method in the case of Lipschitz functions. In Section \ref{sec:theorems} we derive some properties of \ref{opt:T3}. Section \ref{sec:experiments} presents numerical experiments as illustration.

 \section{Scalarization With Competitiveness Method} \label{sec:task_statement}

Let us define the scalarization-type approach  \ref{opt:T1} based on the concept $\gamma$-competitiveness from definition \ref{def:gamma_competitive}. This optimization problem aligns with how problems are typically framed in practice: given historical data on system performance, metrics are calculated. The goal is to find a solution that satisfies the new properties while minimizing the degradation of metric values in previously observed scenarios.

The historical data is represented by the metric values $f_i$ and the solutions $x_i$ for which they were computed. 
Here we denote $v_i := f_i(x_i)$, and, using definition \ref{def:gamma_competitive},  we formulate the optimization problem as follows:



\begin{align*}
    &\min_{x \in K, \quad \gamma \in R_+} \gamma, \quad \quad \tag{SWCM}\label{opt:T1} \\    
    &\text{s.t. } \quad f_i(x) \leq (1 + \gamma) v_i, ~~ \text{for}~~ i = 1, \ldots, m.
\end{align*}

Thus, the goal of the procedure is to find a relative $\gamma$-competitive point that satisfies the required properties while providing the best possible value for the parameter $\gamma$. The formulated problem contains no additional parameters, making it both more specific and interpretable. Informally, a solution with a parameter $\gamma$ can be considered "good enough" based on our current understanding of the objective function values.

\subsection{Related Work}

The approach in \ref{opt:T1} iis closely related to the method proposed in \cite{gembicki1975approach}, which formulates the following optimization problem:

\begin{align*}
    \min_{x, \gamma} \gamma, & \label{\Tilde{opt:T1}} \\    
    \text{s.t. } & x \in K, \\
                 &f_i(x) - w_i \gamma \leq f_i^*, ~~ \text{for}~~ i = 1, \ldots, m.
\end{align*}

Here, \( f_i^* \) are interpreted as target values for the optimized functions, which are not necessarily associated with any particular point and may be chosen based on other considerations. The parameters \( w_i \) represent the relative importance of changes in the \( i \)-th function, essentially determining the direction in which the functions can vary.

This problem can be reduced to our formulation by setting the parameters \( w_i = f_i^* = v_i \). The key distinction between our approach and this one is that we associate the optimization problem with the function values at specific points and fix the parameters. Fixing these parameters enhances the interpretability of the solution. The selected values at given points are then used to approximate a solution for the case of Lipschitz functions, which we now turn to.

\section{\ref{opt:T3}}\label{sec:method} 

\subsection{Lipschitz continuous functions}

Assume that for each \( i = 1, \ldots, m \), the objective function \( f_i(x) : \mathbb{R}^n \rightarrow \mathbb{R}_{++} \) (which is not necessarily convex) is  Lipschitz on $K$ with respect to the norm $\|\cdot\|$ and has a constant \( M_i \). That is, for any $x, y \in K$ the inequality $|f_i(x) - f_i(y)| \leq M_i \|x-y\|$ holds.

Assume that for each \( i = 1, \ldots, m \) the objective function $f_i(\cdot)$ is computed at point \( x_i\) and takes the value \(f_i(x_i) = v_i \).

Now consider some \( x \in K \) and \(\gamma \geq 0\) such that:

\begin{equation*}
    \|x_i - x\| \leq \frac{\gamma v_i}{M_i}.
\end{equation*}

Then we can use the Lipschitz condition to get:
$$
|v_i - f_i(x)| \leq \gamma v_i.
$$

One of the following two alternatives holds::

1. \( v_i \leq f_i(x) \):
    \begin{equation}
        f_i(x) < (1 + \gamma) v_i.
    \end{equation}

2. \( v_i \geq f_i(x) \):
    \begin{equation}
        f_i(x) < v_i < (1 + \gamma) v_i.
    \end{equation}

Thus, the conditions from \ref{opt:T1} are satisfied, and the pair \( x, \gamma \) are feasible points. In this case \( x \) is a relative \(\gamma\)-competitive point for the given functions and points. We can now formulate the approximate optimization problem \ref{opt:T2}.

The solution of this optimization problem is a relative \(\gamma\)-competitive solution for \ref{opt:T0}. Since this is an approximation, the obtained solution may be worse than the exact solution of \ref{opt:T1}.

\begin{align*}
    \min_{x \in K, \gamma\in R_+} \gamma & \tag{$T_2$}\label{opt:T2} \\
    \text{s.t. } &\| x - x_i\| \leq \frac{1}{M_i}(\gamma v_i), ~~ \forall i = 1, \ldots, m.
\end{align*}

We can relax the constraints if additional information about the monotonicity of the functions with respect to their parameters is available. For instance, monotonicity occurs in linear programming problems, where the objective function is monotonic with respect to the parameter \( b \) in the constraint \( Ax \leq b \), as increasing $b$ makes the constraints more lenient. To handle such cases, we define the operator \( \clip_f:\mathbb{R}^n \times \mathbb{R}^n \rightarrow \mathbb{R}^n \) operator for the function \( f \):

\begin{equation}
    \clip_{f}(x, y)_i =         
    \begin{cases}
    \max(x_i - y_i, 0), & f ~\text{increases w.r.t.} ~i\text{-th parameter}, \\
    \max(y_i - x_i, 0), & f ~\text{decreases w.r.t.} ~ i\text{-th parameter}, \\
    x_i - y_i, & \text{otherwise}.
    \end{cases}
\end{equation}

\begin{figure}
\begin{center}    
    \includegraphics[width=0.7\textwidth]{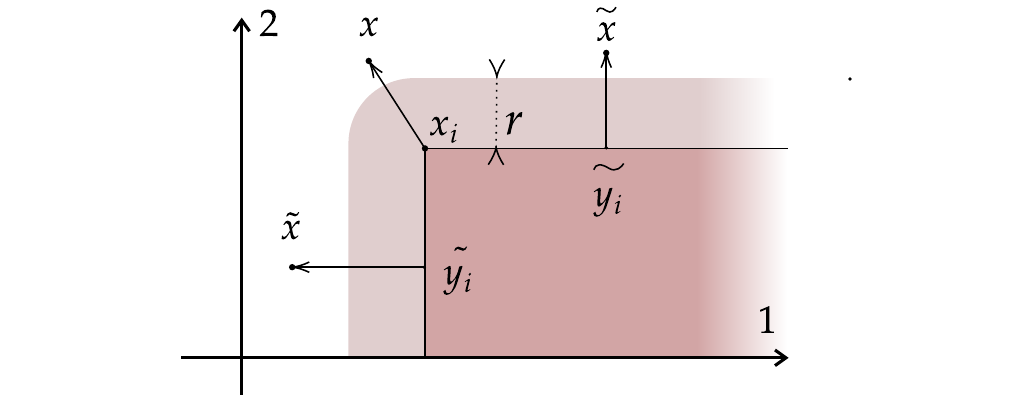}
    \caption{Suppose the function \( f_i \) with reference point $x_i$ decreases with respect to the first parameter and increases with respect to the second parameter, and needs to be minimized. The operator returns \(\mathbf{0}\) for points in the gray area. For other points, the operator returns the vector shown in the figure. The projection is performed coordinately. The shaded area shows the area of space, the projection from which will have a norm not exceeding \( r \), indicated in the figure.}
\end{center}
\label{fig:clip_demo}
\end{figure}

\hyperref[fig:clip_demo]{Figure 1} demonstrates what the operator is doing.
We formulate the optimization problem, in which the difference is replaced by the operator and call it Competitiveness Approximation on Lipschitz Functions (\ref{opt:T3}):

\begin{align*}
    \min_{x, \gamma} \gamma, & \tag{CAoLF}\label{opt:T3} \\
    \text{s.t. } &x \in K, \\
                 &\|\clip_{f_i}(x,x_i)\| \leq \frac{1}{M_i}(\gamma v_i), ~~~ \forall i = 1, \ldots, m.
\end{align*}

\subsection{Functions with Lipschitz Continuous Gradient}

In this subsection, we consider a different type of constraint for functions—convex functions with a Lipschitz continuous gradient. A differentiable function $f$ is said to have an $L$-Lipschitz continuous gradient on a set $K$ if there exists some $L > 0$ such that for all $x, y \in K: \|\nabla f(x) - \nabla f(y)\| \leq L \|x - y\|$. We also recall the definition of convexity: a function $f$ is called convex on $K$ if $\forall x, y \in K: f(x) \geq f(y) + \langle \nabla f(y), x - y \rangle$. A function is said to be concave if the opposite inequality holds. 

To approximate the competitiveness constraint, we can rewrite it in the following form:

\begin{equation}
    f(x) \leq (1 + \gamma) f(x_i) \quad \Rightarrow \quad f(x) - f(x_i) \leq \gamma f(x_i).
\end{equation}

We will constrain the left-hand side of the inequality. For a concave function $f$, we have:

\begin{equation}
f(x) - f(x_i) \leq \langle \nabla f(x_i), x - x_i \rangle \leq \gamma f(x_i).
\end{equation}

If $f$ is convex and has an $L$-Lipschitz continuous gradient, we get:

\begin{equation}
f(x) - f(x_i) \leq \langle \nabla f(x), x - x_i \rangle = \langle \nabla f(x) \pm \nabla f(x_i), x - x_i \rangle.
\end{equation}

This can be bounded as:
\begin{equation}
L \|x - x_i\|^2 + \langle \nabla f(x_i), x - x_i \rangle \leq \gamma f(x_i).
\end{equation}

The resulting constraints take the form of linear and quadratic conditions, requiring only a single evaluation of the gradients and functions at the points $x_i$. These approximated constraints can be used to set up an optimization problem, as in problem \ref{opt:T3}. 

In the following optimization problem, we consider cases where the function $f_i$ is either Lipschitz continuous with constant $M_i$, has an $L$-Lipschitz continuous gradient with constant $L_i$, or is concave. In some instances, multiple constraints may apply to the same function, and these should be preserved since different approximations can have varying levels of accuracy. We call obtained problem \ref{opt:approx}:

\begin{align*}
    \min_{x, \gamma} & \quad \gamma \tag{APPROX} \label{opt:approx} \\
    \text{s.t.} & \quad x \in K, \\
                & \quad \|\clip_{f_i}(x,x_i)\| \leq \frac{1}{M_i} (\gamma v_i), \quad & f_i \text{ is Lipschitz}, \\
                & \quad \langle \nabla f_i(x_i), x - x_i \rangle \leq \gamma v_i, \quad &f_i \text{ is concave}, \\
                & \quad L_i \|x - x_i\|^2 + \langle \nabla f_i(x_i), x - x_i \rangle \leq \gamma v_i, \quad & f_i \text{ is convex with $L$-Lipschitz gradient}, \\
                & \quad \forall i = 1, \dots, m.
\end{align*}

\section{\ref{opt:T3}: Solution properties} \label{sec:theorems}

In this section, we present the theoretical analysis results of the optimization problems discussed. The formulated problems are convex. Theorem 1 demonstrates that the solution of \ref{opt:T3} satisfies the necessary conditions — namely, it is a feasible point for the original problem. Theorem 2 examines the quality of solution achievable when only approximate Lipschitz constants are known, a common scenario in applications due to limited access to function details.

\begin{theorem}[Feasibility]\label{th:feasibility}
    The solution $(x, \gamma)$ of problem \ref{opt:T3} is a feasible point for problem \ref{opt:T1}.
\end{theorem}

{\it Proof}:   
    Let $z_i = \clip_{f_i}(x, x_i)$. Denote $y_i = x - z_i$. This point is a projection of $x$ into gray area in \hyperref[fig:clip_demo]{Figure 1}. For this $y_i$, due to the monotonicity of the function we have $f_i(y_i) \leq f_i(x_i)$. Using the inequality $\| x - y_i \| = \| z_i \| \leq \frac{\gamma v_i}{M_i}$, we obtain $f_i(x) - f_i(y_i) \leq \gamma v_i$. Summing these inequalities gives the desired result: $f_i(x) - f_i(x_i) \leq \gamma v_i$.
\qed

\newcommand{\sto}{\gamma^*} 
\newcommand{\stt}{\widetilde{\gamma}} 
\newcommand{\stf}{\overline{\gamma}} 
\newcommand{\km}{\kappa_{\max}}

\begin{theorem}[Stability]\label{th:stability}
    Let $\widetilde{M_i} = \kappa_i M_i$ denote the approximations of Lipschitz constants for the functions. Denote $\sto$ as the solution to problem \ref{opt:T3}. Let $x$ be the solution to problem \ref{opt:T3_wt}:
    
    \begin{align*}
    \min_{x, \gamma} \gamma, & \tag{$\widetilde{T_3}$}\label{opt:T3_wt} \\
    \text{s.t. } & x \in K, \\
                 & \| \clip_{f_i}(x, x_i) \| \leq \frac{1}{\widetilde{M_i}} (\gamma v_i), ~~~ \forall i = 1, \ldots, m.
    \end{align*}
    For the obtained solution $x$, it holds that:
    $$|f(x) - f(x_i)| \leq \frac{\kappa_{\max}}{\kappa_i} \sto v_i$$
    where $\kappa_{\max} = \max_{i = 1, \ldots, m} \kappa_i$.
\end{theorem}

\begin{remark}
    If all constants are multiplied by the same factor, we obtain the same $x$ as for the problem with exact constants. Significant deterioration in function value may occur if some constants are estimated much worse than others. 
\end{remark}
\begin{remark}
    Problem \ref{opt:T3_wt} differs from \ref{opt:T3} only by replacing Lipschitz constants with approximations.
\end{remark}

{\it Proof} 
    Let's introduce an optimization problem where all Lipschitz constants are multiplied by $\kappa_{\max}$:
    
    \begin{align*}
    \min_{x, \gamma} \gamma, & \tag{$T_4$}\label{opt:T4} \\
    \text{s.t. } & x \in K, \\
                 & \| \clip_{f_i}(x, x_i) \| \leq \frac{1}{\kappa_{\max} M_i} (\gamma v_i), ~~~ \forall i = 1, \ldots, m.
    \end{align*}    
    Denote $\stt$ as the solution to \ref{opt:T3} with constants $\widetilde{L_i}$, and $\stf$ as the solution to problem \ref{opt:T4}.
    
    Note that if we substitute $\gamma = a$ and $\frac{\gamma}{\kappa_{\max}} = a$ in problems \ref{opt:T3} and \ref{opt:T4} respectively, we obtain the same problem. Hence, $\sto = \frac{\stf}{\kappa_{\max}}$.
    
    For problems $\widetilde{T_3}$ and \ref{opt:T4}, the relation $\stt \leq \stf$ holds: we can consider the intersection of the spheres in the solution of \ref{opt:T4} with parameter $\stf$. If we increase the radii of the spheres, the intersection becomes larger and we can reduce $\gamma$, which happens in $\widetilde{T_3}$.
    
    Thus, we get $\stt \leq \kappa_{\max} (\sto)$ and from the conditions on the radii in $\widetilde{T_3}$, $x$ satisfies:
    \begin{equation}
        |f(x) - f(x_i)| \leq \frac{M_i}{M_i \kappa_i} (\stt) \leq \frac{\kappa_{\max}}{\kappa_i} \sto.
    \end{equation}
\qed

\begin{remark}[Optimality of the obtained estimate]
    To evaluate the quality degradation of the approximate solution, consider an example with two functions: $f_1(x) = 1 + M_1 |x|$ and $f_2(x) = 1 + M_2 |x - 1|$ with $x_1 = 0, x_2 = 1$.
    
    For \ref{opt:T3}, we obtain the solution $\sto: (\sto) (\frac{1}{M_1} + \frac{1}{M_2}) = 1$. For the same problem with approximate constants $\widetilde{M_i} = \kappa_i M_i$, we get $\stt: (\stt) (\frac{1}{M_1 \kappa_1} + \frac{1}{M_2 \kappa_2}) = 1$.
    
    Hence, for the case where $M_2 << M_1$ and $\kappa_1 = 1$ (know exact $M_1$):
    \begin{equation}
        \frac{\stt}{\sto} = 
        \frac{\frac{1}{M_1} + \frac{1}{M_2}}{\frac{1}{M_1 \kappa_1} + \frac{1}{M_2 \kappa_2}} = 
        \frac{\frac{M_2}{M_1} + 1}{\frac{M_2}{M_1 \kappa_1} + \frac{1}{\kappa_2}} \approx \kappa_2.
    \end{equation}
    Thus, the obtained solution will be close to the boundary obtained in the theorem. Figure 2 visually illustrates the deviation of the point found by the algorithm due to incorrect estimation of constants.
\end{remark}

\begin{figure}[ht]
\begin{center}    
        \includegraphics[width = 0.7\textwidth]{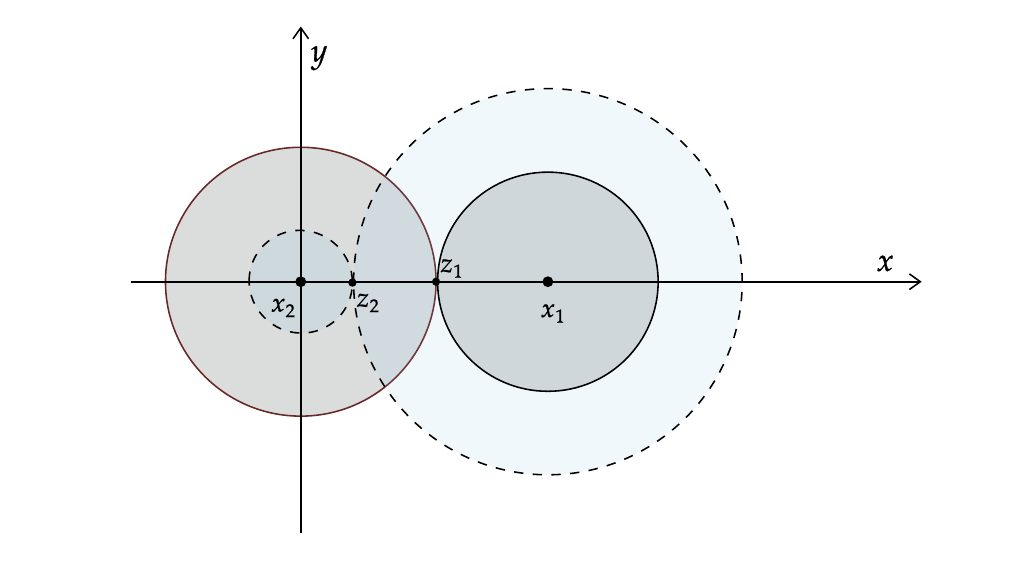}
    \caption{Result of approximate Lipschitz constant estimation. Intersection of gray circles $z_1$ — exact solution. However, the algorithm finds $z_2$ due to incorrect constant estimation.}
\end{center}
\label{fig:theorem_example}
\end{figure}

\section{Numerical Experiments} \label{sec:experiments}
To demonstrate the work of \ref{opt:T3}, two experiments were conducted. The first experiment deals with the Minimum Cost Concurrent Flow (MCCF) problem. In this experiment, the objective functions have a similar structure but differ in parameters corresponding to network usage scenarios. The second experiment consider different metrics: MCCF, Maximum Flow and $\lambda_2$ value of graph. The code for the experiments can be found at \url{https://anonymous.4open.science/r/paper1-A2DF}.

\subsection{Experiment with MCCF function }\label{exp:concurrent}

MCCF is stated as follows: given a graph $\mathcal{G}(V, E)$ with $k$ vertices and $n$ edges, where each edge $e \in E$ has a capacity $b_e$ and a cost per unit of flow $c_e$. Additionally, there are flows that need to be managed through the network: $(s_i, t_i, f_i)$. This tuple represents a source $s$, a target $t$, and the required flow amount $f$. These requests are represented by a demand matrix $D$: $\forall i: D_{s_i, t_i} = f_i$. The goal is to find the flows in the graph $F \in\mbR^{n\times k}$ that satisfy all demands at the minimum possible cost.

The problem is reformulated as a linear programming problem \cite{bazaraa2011linear}. We define the vertex-edge incidence matrix $A \in \mbR^{k\times n}$. Each column contains exactly two non-zero values. The column corresponding to edge $(i, j)$ has "+1" in row $i$, "-1" in row $j$, and zeros in all other rows. For each vertex, a supply vector $d_i \in \mbR^k$ is introduced to record all outgoing flows from this vertex. Using the demand matrix $D$, the components of $d_i$ are expressed as follows:
\begin{equation*}
d_{ij} = 
\begin{cases}
    \sum_{k \neq j} D_{ik}, & i = j, \\
    -D_{ij}, & i \neq j.
\end{cases}
\end{equation*}

We combine these column vectors into a matrix denoted as $D_A \in \mbR^{k\times k}$, which is used for convenient notation of the problem. Using \cite{bazaraa2011linear}, we formulate the optimization problem as follows:

\begin{align*}
    \min_{F} c_e^T F \textbf{1}, & \\
    \text{s.t.} & ~ F\textbf{1} \leq b, \\
                & ~ AF = D_A.\\
\end{align*}

In this formulation, the problem may not have a feasible solution due to the constraint with $b$. To remove the constraints on $b$ and simplify the problem, we add a penalty $y$:

\begin{align*}
    \min_{F} c^T F \textbf{1} + c_a^T y, & \tag{$MCCF$}\label{opt:MCCF}\\
    \text{s.t.} & ~ F\textbf{1} \leq b + y,\\
                & ~ AF = D_A.\\
\end{align*}

This penalty is interpreted as renting additional capacity in the network. We denote $f_D(b) = \text{value}(MCCF(b, D))$,.

Now let's describe a scenario of using the method. There is a company that provides delivery services in a certain network. It needs to rent bandwidth for its operations. The company can rent bandwidth at the beginning of a period at a price $c_b$ and rent additional bandwidth during the period at a price $c_a$. In practice, this occurs in stages since the bandwidth is rented in advance for a long term. At each time period, there are supply demands — the demand matrix $D_i$, which must be fulfilled. These matrices define the network usage scenarios. The company is unaware of these demands at the start of the period. Suppose it operated for $m$ periods with different rented capacities $b_i$ and observed $D_i$. At the end of the period, it evaluates its expenses during the period $f_i(b_i) = f_{D_i}(b_i)$. These functions satisfy Lipschitz conditions. The goal is to purchase capacity within a given budget to ensure minimal costs for different scenarios. The company only knows the $D_i$ it has observed so far. This necessitates solving problem \ref{opt:T3}. Constraints are defined as $K = \{b: c_b^T b \leq B\}$.

\subsubsection{Numerical results}

In the experiment, we use the "germany50" topology from the SNDLib dataset \cite{orlowski2010sndlib}. This is a network of 50 vertices and 176 edges. The dataset includes a set of demand matrices and edge flow costs. These matrices will be used as $D$ matrices. To ensure problem feasibility with smaller capacities, we sparse the demands: with a probability of 0.4, an entry in the matrix is zeroed out. Parameters for the experiments were generated as follows: flow costs $c$ are given in the graphs. Based on these, initial rental costs and in-period rental costs were generated. If a capacity is cheaper, it is a higher quality capacity, so its cost should be higher. Also, we assume renting in-period is more expensive than pre-period renting. Thus, pre-period costs were generated as $(c_b)_i = (C/\sqrt{c_i}) \xi$, where $\xi \sim U[9, 11]$ and C is some constant. In-period rental costs are higher: $(c_a)_i = (c_b)_i * \xi, \sim U[1.05, 1.15]$, i.e., 5-15\% more expensive.

\begin{figure}
\centering
\includegraphics[width=\textwidth]{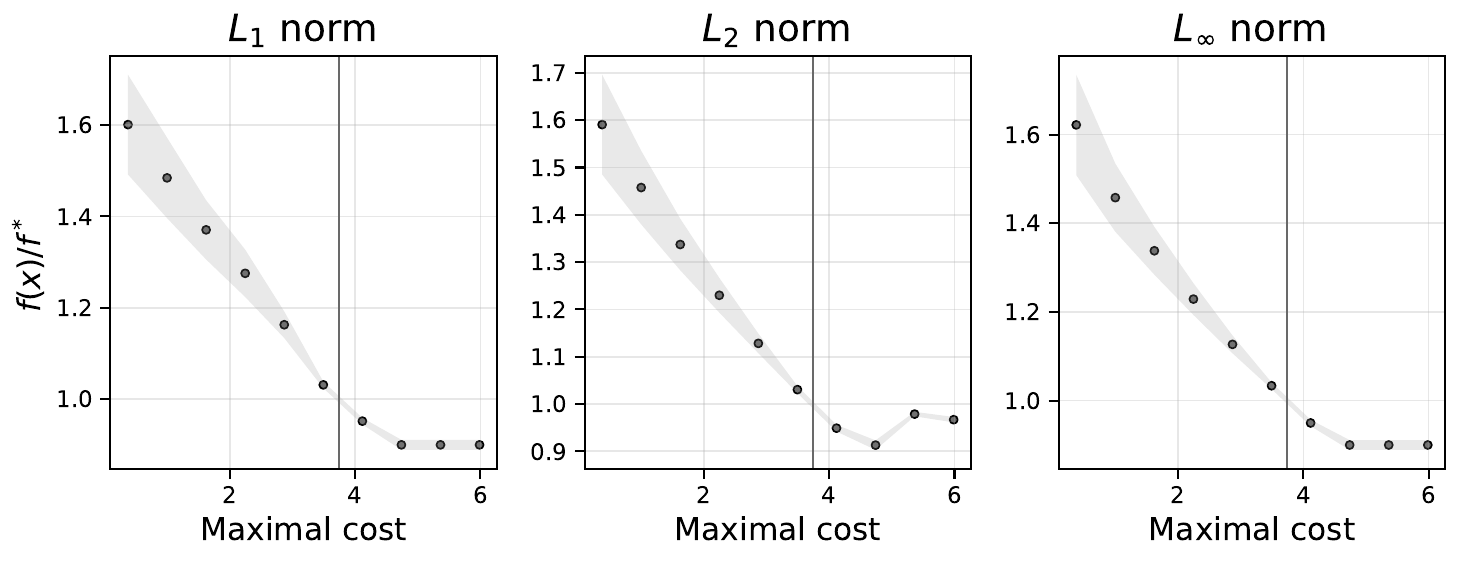}
\caption{The dependence of the relative increase in the value of the function for solving the problem \ref{opt:T3} with MCCF functions on the budget for various norms.}
\label{fig:mcf:relations}
\end{figure}

For this problem, we could not find an exact solution as the solver did not converge. We considered norms $\|\cdot\|_1$; $\|\cdot\|_2$; $\|\cdot\|_\infty$. Figure \ref{fig:mcf:relations} shows the relative increase in function values for different budget values. The point is the average ratio of function values to initial point values for considered budgets. The vertical line represents the average value of the initial budget in the scenarios considered. We tested the algorithm at 10 different budget values from $0.1$ to $1.6$ times of average initial budget. The shaded area is $\pm$ variance. 

\begin{figure}
    \centering
    \includegraphics[width=0.5\textwidth]{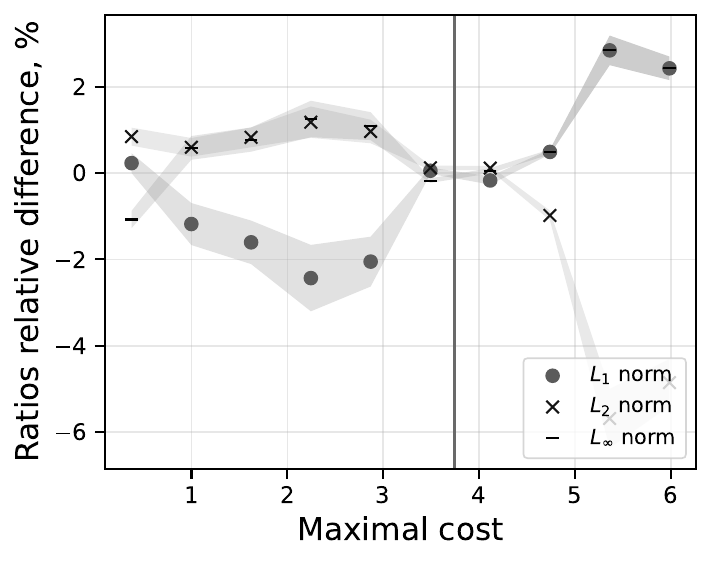}
    \caption{The relative difference in quality of the approximate solution compared to the average quality solution for the solution of the problem \ref{opt:T3} with MCCF functions. We see that different norms work better on different budgets. For example, the $L_1$ norm works better on large budgets because it buys less extra resources.}
    \label{fig:mcf:relative_diff}
\end{figure}

Since there is no exact solution for this case, we compare these solutions with the average quality of solutions. This is shown in Figure \ref{fig:mcf:relative_diff}. Higher values are better. We see that different norms work better for different budgets: the $L_1$ norm gives sparse solutions, which are not optimal for small budgets as some resources are minimally purchased but are good for large budgets as excess resources are the cheapest. The $L_2$ norm yields a solution that is better than the average quality of solution for small budgets. Similarly to $L_1$, but already for larger budgets.

\subsection{Experiment with different functions}

In this section, we consider different functions to be optimized. One of these functions is the MCCF, described above. The second function is the maximum flow function, which is important for ensuring the ability to send large amounts of information within the network. The last metric is the $\lambda_2$ value of the graph's Laplacian. Works, such as \cite{mohar1991laplacian}, demonstrate that the value of $\lambda_2$ is connected with various graph invariants, including connectivity, expanding properties, the isoperimetric number, the maximum cut, the independence number, genus, diameter, mean distance, and bandwidth-type parameters of a graph.

\subsubsection{Maximum flow}
Let's consider the maximum flow function. It is formulated as a linear optimization problem \cite{bazaraa2011linear}:
\begin{align*}
    \max_{f, x} f, &\\
    \text{s.t.}~~ &(e_t - e_s)f + Ax = 0, \\
    &x \leq b, \\
    &x \geq 0.
\end{align*}

Here, $f\in\mbR$ represents the amount of flow that can be sent. $x\in\mbR^n$ represents the flows within the edges. They must be non-negative and are constrained by bandwidth. $e_s$ and $e_t$ are one-hot vectors with ones in the appropriate positions for the source and target, respectively.

\subsubsection{ \texorpdfstring{$\lambda_2$}{Lg} value}

$\lambda_2$ is defined for undirected graphs. Let $W$ be the weighted adjacency matrix, which is symmetric. The value of $w_{ij}$ is the bandwidth of the edge $(i,j)$ if there is an edge, and 0 otherwise. The Laplacian matrix of a graph is defined as $L_G = diag(W\textbf{1}) - W$. The eigenvalues of $L_G$ are non-negative, and the smallest eigenvalue is 0, corresponding to the vector $\textbf{1}$. The second smallest eigenvalue of the Laplacian is $\lambda_2$. It is related to various characteristics of the graph. Higher $\lambda_2$ values indicate better graph properties, so our goal is to maximize this value. More details about the properties of the Laplacian and $\lambda_2$ can be found in \cite{mohar1991laplacian}.

\subsubsection{Numerical results}

The data generation process and tested budgets selection is the same as in the previous experiment. For the maximal flow function, we need to define source-target pairs. This is done as follows: for each demand matrix realization, we take 5\% of the largest queries and calculate the maximal flow for these pairs. This approach is natural since intensive exchanges are expected for these pairs in subsequent time periods.

\begin{figure}
    \centering
    \includegraphics[width=0.9\textwidth]{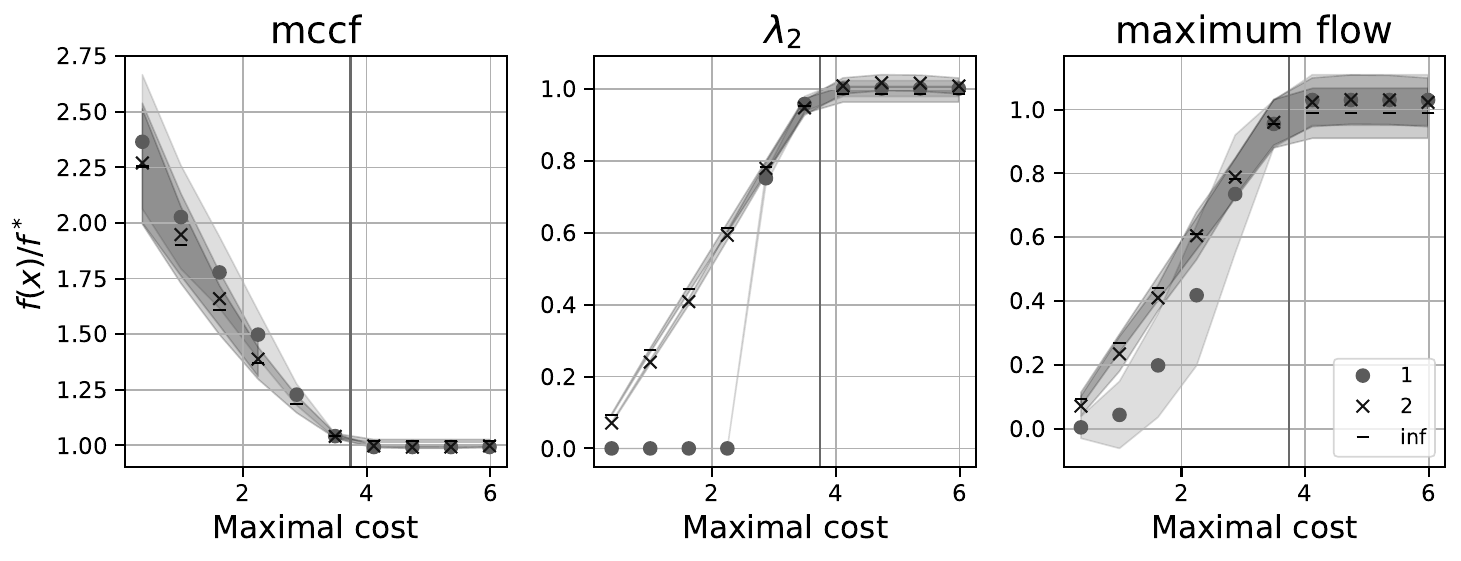}
    \caption{The dependence of the relative increase/decrease of considered functions for solutions of \ref{opt:T3} w.r.t. different norms. The vertical line indicates the mean budget used for $b_i$ in previous periods. We see that the $L_1$ norm does not work well for the $\lambda_2$ metric with small budgets. The $L_2$ and $L_\infty$ norms work comparably on all the functions considered.}
    \label{fig:diff:relations}
\end{figure}

Figure \ref{fig:diff:relations} shows the results of solving \ref{opt:T3} for considered functions with different norms. The x-axis represents the budgets. The vertical line indicates the mean budget used for $b_i$ in previous periods. The y-axis shows the ratio of function values at the resulting solution for a given budget to the function values at points $b_i$. The dot represents the mean, and the colored area represents $\pm$ variance. The MCCF function is minimized, so its graph is decreasing, while the other functions are maximized, so their graphs are increasing.

The results show that the $L_1$ norm performs the worst for budgets lower than the mean used budgets. However, for large budgets, it performs better than the $L_\infty$ norm. Among the considered norms, $L_2$  norm shows the best results.

\section{Conclusions}
This work proposes a scalarization method \ref{opt:T1} for multi-criteria optimization which does not require hyperparameter tuning. The method is based on a generalization of the definition of competitive solutions. We also established a connection between our scalarization method and the approach proposed in \cite{gembicki1975approach}. Our method is a special case of this approach and, due to the fixed parameters, its solutions have a clear interpretation as relative $\gamma$-competitive solutions.

For Lipschitz and Lipschitz monotonic functions, we presented \ref{opt:T3} -- optimization problem to search approximate solution of \ref{opt:T1}. This method is particularly useful when function evaluations are expensive, as it operates solely with a precomputed set of values. The necessity for such methods is driven by the scale of modern problems, such as in the optimization of TV network topologies and supply chains. The conducted experiments demonstrated the effectiveness and practical utility of the proposed method, confirming its applicability in various scenarios.

\bibliography{main}


\begin{thebibliography}{24}
\ifx \bisbn   \undefined \def \bisbn  #1{ISBN #1}\fi
\ifx \binits  \undefined \def \binits#1{#1}\fi
\ifx \bauthor  \undefined \def \bauthor#1{#1}\fi
\ifx \batitle  \undefined \def \batitle#1{#1}\fi
\ifx \bjtitle  \undefined \def \bjtitle#1{#1}\fi
\ifx \bvolume  \undefined \def \bvolume#1{\textbf{#1}}\fi
\ifx \byear  \undefined \def \byear#1{#1}\fi
\ifx \bissue  \undefined \def \bissue#1{#1}\fi
\ifx \bfpage  \undefined \def \bfpage#1{#1}\fi
\ifx \blpage  \undefined \def \blpage #1{#1}\fi
\ifx \burl  \undefined \def \burl#1{\textsf{#1}}\fi
\ifx \doiurl  \undefined \def \doiurl#1{\url{https://doi.org/#1}}\fi
\ifx \betal  \undefined \def \betal{\textit{et al.}}\fi
\ifx \binstitute  \undefined \def \binstitute#1{#1}\fi
\ifx \binstitutionaled  \undefined \def \binstitutionaled#1{#1}\fi
\ifx \bctitle  \undefined \def \bctitle#1{#1}\fi
\ifx \beditor  \undefined \def \beditor#1{#1}\fi
\ifx \bpublisher  \undefined \def \bpublisher#1{#1}\fi
\ifx \bbtitle  \undefined \def \bbtitle#1{#1}\fi
\ifx \bedition  \undefined \def \bedition#1{#1}\fi
\ifx \bseriesno  \undefined \def \bseriesno#1{#1}\fi
\ifx \blocation  \undefined \def \blocation#1{#1}\fi
\ifx \bsertitle  \undefined \def \bsertitle#1{#1}\fi
\ifx \bsnm \undefined \def \bsnm#1{#1}\fi
\ifx \bsuffix \undefined \def \bsuffix#1{#1}\fi
\ifx \bparticle \undefined \def \bparticle#1{#1}\fi
\ifx \barticle \undefined \def \barticle#1{#1}\fi
\bibcommenthead
\ifx \bconfdate \undefined \def \bconfdate #1{#1}\fi
\ifx \botherref \undefined \def \botherref #1{#1}\fi
\ifx \url \undefined \def \url#1{\textsf{#1}}\fi
\ifx \bchapter \undefined \def \bchapter#1{#1}\fi
\ifx \bbook \undefined \def \bbook#1{#1}\fi
\ifx \bcomment \undefined \def \bcomment#1{#1}\fi
\ifx \oauthor \undefined \def \oauthor#1{#1}\fi
\ifx \citeauthoryear \undefined \def \citeauthoryear#1{#1}\fi
\ifx \endbibitem  \undefined \def \endbibitem {}\fi
\ifx \bconflocation  \undefined \def \bconflocation#1{#1}\fi
\ifx \arxivurl  \undefined \def \arxivurl#1{\textsf{#1}}\fi
\csname PreBibitemsHook\endcsname

\bibitem[\protect\citeauthoryear{Miettinen}{1999}]{miettinen1999nonlinear}
\begin{bbook}
\bauthor{\bsnm{Miettinen}, \binits{K.}}:
\bbtitle{Nonlinear Multiobjective Optimization}
vol. \bseriesno{12}.
\bpublisher{Springer},
\blocation{New York}
(\byear{1999})
\end{bbook}
\endbibitem

\bibitem[\protect\citeauthoryear{Greco}{2006}]{greco2006multiple}
\begin{bbook}
\bauthor{\bsnm{Greco}, \binits{S.}}:
\bbtitle{Multiple Criteria Decision Analysis: State of the Art Surveys}.
\bsertitle{International Series in Operations Research \& Management Science}.
\bpublisher{Springer},
\blocation{New York}
(\byear{2006}).
\burl{https://books.google.ru/books?id=lzyNasqOxAQC}
\end{bbook}
\endbibitem

\bibitem[\protect\citeauthoryear{K{\"o}ksalan et~al.}{2011}]{koksalan2011multiple}
\begin{bchapter}
\bauthor{\bsnm{K{\"o}ksalan}, \binits{M.}},
\bauthor{\bsnm{Wallenius}, \binits{J.}},
\bauthor{\bsnm{Zionts}, \binits{S.}}:
\bctitle{Multiple criteria decision making: From early history to the 21st century}.
(\byear{2011}).
\burl{https://api.semanticscholar.org/CorpusID:109047337}
\end{bchapter}
\endbibitem

\bibitem[\protect\citeauthoryear{Altiparmak et~al.}{2006}]{altiparmak2006genetic}
\begin{barticle}
\bauthor{\bsnm{Altiparmak}, \binits{F.}},
\bauthor{\bsnm{Gen}, \binits{M.}},
\bauthor{\bsnm{Lin}, \binits{L.}},
\bauthor{\bsnm{Paksoy}, \binits{T.}}:
\batitle{A genetic algorithm approach for multi-objective optimization of supply chain networks}.
\bjtitle{Computers \& industrial engineering}
\bvolume{51}(\bissue{1}),
\bfpage{196}--\blpage{215}
(\byear{2006})
\end{barticle}
\endbibitem

\bibitem[\protect\citeauthoryear{Elmusrati et~al.}{2008}]{elmusrati2008applications}
\begin{barticle}
\bauthor{\bsnm{Elmusrati}, \binits{M.}},
\bauthor{\bsnm{El-Sallabi}, \binits{H.}},
\bauthor{\bsnm{Koivo}, \binits{H.}}:
\batitle{Applications of multi-objective optimization techniques in radio resource scheduling of cellular communication systems}.
\bjtitle{IEEE Transactions on Wireless Communications}
\bvolume{7}(\bissue{1}),
\bfpage{343}--\blpage{353}
(\byear{2008})
\end{barticle}
\endbibitem

\bibitem[\protect\citeauthoryear{Mastrocinque et~al.}{2013}]{mastrocinque2013multi}
\begin{barticle}
\bauthor{\bsnm{Mastrocinque}, \binits{E.}},
\bauthor{\bsnm{Yuce}, \binits{B.}},
\bauthor{\bsnm{Lambiase}, \binits{A.}},
\bauthor{\bsnm{Packianather}, \binits{M.S.}}:
\batitle{A multi-objective optimization for supply chain network using the bees algorithm}.
\bjtitle{International Journal of Engineering Business Management}
\bvolume{5},
\bfpage{38}
(\byear{2013})
\end{barticle}
\endbibitem

\bibitem[\protect\citeauthoryear{Bjornson et~al.}{2014}]{bjornson2014multiobjective}
\begin{barticle}
\bauthor{\bsnm{Bjornson}, \binits{E.}},
\bauthor{\bsnm{Jorswieck}, \binits{E.A.}},
\bauthor{\bsnm{Debbah}, \binits{M.}},
\bauthor{\bsnm{Ottersten}, \binits{B.}}:
\batitle{Multiobjective signal processing optimization: The way to balance conflicting metrics in 5g systems}.
\bjtitle{IEEE Signal Processing Magazine}
\bvolume{31}(\bissue{6}),
\bfpage{14}--\blpage{23}
(\byear{2014})
\end{barticle}
\endbibitem

\bibitem[\protect\citeauthoryear{Suttorp and Igel}{2006}]{suttorp2006multi}
\begin{botherref}
\oauthor{\bsnm{Suttorp}, \binits{T.}},
\oauthor{\bsnm{Igel}, \binits{C.}}:
Multi-objective optimization of support vector machines.
Multi-objective machine learning,
199--220
(2006)
\end{botherref}
\endbibitem

\bibitem[\protect\citeauthoryear{Zuluaga et~al.}{2013}]{zuluaga2013active}
\begin{bchapter}
\bauthor{\bsnm{Zuluaga}, \binits{M.}},
\bauthor{\bsnm{Sergent}, \binits{G.}},
\bauthor{\bsnm{Krause}, \binits{A.}},
\bauthor{\bsnm{P{\"u}schel}, \binits{M.}}:
\bctitle{Active learning for multi-objective optimization}.
In: \bbtitle{International Conference on Machine Learning},
pp. \bfpage{462}--\blpage{470}
(\byear{2013}).
\bcomment{PMLR}
\end{bchapter}
\endbibitem

\bibitem[\protect\citeauthoryear{Sener and Koltun}{2018}]{sener2018multi}
\begin{botherref}
\oauthor{\bsnm{Sener}, \binits{O.}},
\oauthor{\bsnm{Koltun}, \binits{V.}}:
Multi-task learning as multi-objective optimization.
Advances in neural information processing systems
\textbf{31}
(2018)
\end{botherref}
\endbibitem

\bibitem[\protect\citeauthoryear{Rangaiah and Petriciolet}{2013}]{rangaiah2013multi}
\begin{botherref}
\oauthor{\bsnm{Rangaiah}, \binits{G.P.}},
\oauthor{\bsnm{Petriciolet}, \binits{A.}}:
Multi-objective optimization in chemical engineering.
Developments and applications/edited by Gade Pandu Rangaiah, Adri{\'a}n Bonilla-Petriciolet
(2013)
\end{botherref}
\endbibitem

\bibitem[\protect\citeauthoryear{Boada et~al.}{2016}]{boada2016multi}
\begin{barticle}
\bauthor{\bsnm{Boada}, \binits{Y.}},
\bauthor{\bsnm{Reynoso-Meza}, \binits{G.}},
\bauthor{\bsnm{Pic{\'o}}, \binits{J.}},
\bauthor{\bsnm{Vignoni}, \binits{A.}}:
\batitle{Multi-objective optimization framework to obtain model-based guidelines for tuning biological synthetic devices: an adaptive network case}.
\bjtitle{BMC systems biology}
\bvolume{10},
\bfpage{1}--\blpage{19}
(\byear{2016})
\end{barticle}
\endbibitem

\bibitem[\protect\citeauthoryear{Marler and Arora}{2004}]{marler2004survey}
\begin{barticle}
\bauthor{\bsnm{Marler}, \binits{R.T.}},
\bauthor{\bsnm{Arora}, \binits{J.S.}}:
\batitle{Survey of multi-objective optimization methods for engineering}.
\bjtitle{Structural and multidisciplinary optimization}
\bvolume{26},
\bfpage{369}--\blpage{395}
(\byear{2004})
\end{barticle}
\endbibitem

\bibitem[\protect\citeauthoryear{Ngatchou et~al.}{2005}]{ngatchou2005pareto}
\begin{bchapter}
\bauthor{\bsnm{Ngatchou}, \binits{P.}},
\bauthor{\bsnm{Zarei}, \binits{A.}},
\bauthor{\bsnm{El-Sharkawi}, \binits{A.}}:
\bctitle{Pareto multi objective optimization}.
In: \bbtitle{Proceedings of the 13th International Conference On, Intelligent Systems Application to Power Systems},
pp. \bfpage{84}--\blpage{91}
(\byear{2005}).
\bcomment{IEEE}
\end{bchapter}
\endbibitem

\bibitem[\protect\citeauthoryear{Konak et~al.}{2006}]{konak2006multi}
\begin{barticle}
\bauthor{\bsnm{Konak}, \binits{A.}},
\bauthor{\bsnm{Coit}, \binits{D.W.}},
\bauthor{\bsnm{Smith}, \binits{A.E.}}:
\batitle{Multi-objective optimization using genetic algorithms: A tutorial}.
\bjtitle{Reliability engineering \& system safety}
\bvolume{91}(\bissue{9}),
\bfpage{992}--\blpage{1007}
(\byear{2006})
\end{barticle}
\endbibitem

\bibitem[\protect\citeauthoryear{Miettinen et~al.}{2008}]{miettinen2008introduction}
\begin{bchapter}
\bauthor{\bsnm{Miettinen}, \binits{K.}},
\bauthor{\bsnm{Ruiz}, \binits{F.}},
\bauthor{\bsnm{Wierzbicki}, \binits{A.P.}}:
\bctitle{Introduction to multiobjective optimization: interactive approaches}.
In: \bbtitle{Multiobjective Optimization: Interactive and Evolutionary Approaches},
pp. \bfpage{27}--\blpage{57}.
\bpublisher{Springer},
\blocation{New York}
(\byear{2008})
\end{bchapter}
\endbibitem

\bibitem[\protect\citeauthoryear{Daulton et~al.}{2022}]{daulton2022multi}
\begin{bchapter}
\bauthor{\bsnm{Daulton}, \binits{S.}},
\bauthor{\bsnm{Eriksson}, \binits{D.}},
\bauthor{\bsnm{Balandat}, \binits{M.}},
\bauthor{\bsnm{Bakshy}, \binits{E.}}:
\bctitle{Multi-objective bayesian optimization over high-dimensional search spaces}.
In: \bbtitle{Uncertainty in Artificial Intelligence},
pp. \bfpage{507}--\blpage{517}
(\byear{2022}).
\bcomment{PMLR}
\end{bchapter}
\endbibitem

\bibitem[\protect\citeauthoryear{Borodin and El-Yaniv}{1998}]{borodin2005online}
\begin{bchapter}
\bauthor{\bsnm{Borodin}, \binits{A.}},
\bauthor{\bsnm{El-Yaniv}, \binits{R.}}:
\bctitle{Online computation and competitive analysis}.
(\byear{1998}).
\burl{https://api.semanticscholar.org/CorpusID:5431684}
\end{bchapter}
\endbibitem

\bibitem[\protect\citeauthoryear{Banos and Papageorgiou}{1995}]{banos1995linear}
\begin{barticle}
\bauthor{\bsnm{Banos}, \binits{J.M.}},
\bauthor{\bsnm{Papageorgiou}, \binits{M.}}:
\batitle{A linear programming approach to large-scale linear optimal control problems}.
\bjtitle{IEEE transactions on automatic control}
\bvolume{40}(\bissue{5}),
\bfpage{971}--\blpage{977}
(\byear{1995})
\end{barticle}
\endbibitem

\bibitem[\protect\citeauthoryear{Martin}{2012}]{martin2012large}
\begin{bbook}
\bauthor{\bsnm{Martin}, \binits{R.K.}}:
\bbtitle{Large Scale Linear and Integer Optimization: a Unified Approach}.
\bpublisher{Springer},
\blocation{New York}
(\byear{2012})
\end{bbook}
\endbibitem

\bibitem[\protect\citeauthoryear{Gembicki and Haimes}{1975}]{gembicki1975approach}
\begin{barticle}
\bauthor{\bsnm{Gembicki}, \binits{F.}},
\bauthor{\bsnm{Haimes}, \binits{Y.}}:
\batitle{Approach to performance and sensitivity multiobjective optimization: The goal attainment method}.
\bjtitle{IEEE Transactions on Automatic control}
\bvolume{20}(\bissue{6}),
\bfpage{769}--\blpage{771}
(\byear{1975})
\end{barticle}
\endbibitem

\bibitem[\protect\citeauthoryear{Bazaraa et~al.}{2011}]{bazaraa2011linear}
\begin{bbook}
\bauthor{\bsnm{Bazaraa}, \binits{M.S.}},
\bauthor{\bsnm{Jarvis}, \binits{J.J.}},
\bauthor{\bsnm{Sherali}, \binits{H.D.}}:
\bbtitle{Linear Programming and Network Flows}.
\bpublisher{John Wiley \& Sons},
\blocation{Hoboken, New Jersey}
(\byear{2011})
\end{bbook}
\endbibitem

\bibitem[\protect\citeauthoryear{Orlowski et~al.}{2010}]{orlowski2010sndlib}
\begin{barticle}
\bauthor{\bsnm{Orlowski}, \binits{S.}},
\bauthor{\bsnm{Wess{\"a}ly}, \binits{R.}},
\bauthor{\bsnm{Pi{\'o}ro}, \binits{M.}},
\bauthor{\bsnm{Tomaszewski}, \binits{A.}}:
\batitle{Sndlib 1.0—survivable network design library}.
\bjtitle{Networks: An International Journal}
\bvolume{55}(\bissue{3}),
\bfpage{276}--\blpage{286}
(\byear{2010})
\end{barticle}
\endbibitem

\bibitem[\protect\citeauthoryear{Mohar et~al.}{1991}]{mohar1991laplacian}
\begin{barticle}
\bauthor{\bsnm{Mohar}, \binits{B.}},
\bauthor{\bsnm{Alavi}, \binits{Y.}},
\bauthor{\bsnm{Chartrand}, \binits{G.}},
\bauthor{\bsnm{Oellermann}, \binits{O.}}:
\batitle{The laplacian spectrum of graphs}.
\bjtitle{Graph theory, combinatorics, and applications}
\bvolume{2}(\bissue{871-898}),
\bfpage{12}
(\byear{1991})
\end{barticle}
\endbibitem

\end{thebibliography}

\end{document}